\documentclass{article}
\usepackage{amssymb,amsmath}
\usepackage{graphics}

\newcommand\acl{\hbox{\rm acl}}

\newcommand{\proof}{\noindent {\bf Proof.}\hspace{2mm}}
\newcommand{\proofof}[1] {\noindent {\bf Proof of #1.}\hspace
{2mm}}

\newcommand{\qed}{{$\square$}\medbreak}

\newenvironment{repstat}[1]{
\noindent
{\bf #1 }\em }

\newtheorem{lemma}{Lemma}[section]
\newtheorem{theorem}{Theorem}
\newtheorem{prop}[lemma]{Proposition}
\newtheorem{fact}[lemma]{Fact}

\newtheorem{cor}[lemma]{Corollary}

\newtheorem{hypothesis}[lemma]{Hypothesis}

\newcommand\normal{\lhd}

\newcommand\sd{\rtimes}
\newcommand\rk{\hbox{\rm rk}\,}

\newcommand\pr{\hbox{\em pr}\,}

\newcommand\GL{\hbox{\rm GL}}

\newcommand\End{\hbox{\rm End}}

\newcommand\Stab{\hbox{\rm Stab}}

\newcommand*\invlim{\varprojlim}

\newfont{\Bbbb}{msbm10}
\newcommand\nn{\mbox {\Bbbb N}}
\newcommand\zz{\mbox {\Bbbb Z}}
\newcommand\qq{\mbox {\Bbbb Q}}

\newcommand\cc{\mbox {\Bbbb C}}

\input amssym.def

\newsymbol\square 1003
\newsymbol\rtimes 226F

\begin{document}

\title{On analogies between algebraic groups and groups
of finite Morley rank}

\author{Tuna Alt\i nel,
Jeffrey Burdges
\thanks{Supported by NSF postdoctoral fellowship DMS-0503036
 and a Bourse Chateaubriand postdoctoral fellowship
and partially by the ANR project ``Groupes, G\'eom\'etrie et Logique''
JC0$5\_4$7037:jaligot:eric:}\\
Universit\'e de Lyon\\ Universit\'e Lyon 1\\
CNRS UMR 5208
Institut Camille Jordan\\
43 blvd du 11 novembre 1918\\
69622 Villeurbanne cedex, France
}
\maketitle

\begin{abstract}
{We prove that in a connected group of finite Morley
rank the centralizers of decent tori are connected. We then
apply this result to the analysis of minimal connected simple
groups of finite Morley rank. Our applications include
general covering properties by Borel subgroups, the description of
Weyl groups and the analysis of toral automorphisms.}
\end{abstract}

\section{Introduction}
Linear algebraic groups and infinite groups
of finite Morley rank have many analogous
properties. This is expressed in the
strongest possible way by the
central problem in the analysis of groups of finite
Morley rank, namely the {\em Cherlin-Zil'ber
algebraicity conjecture} which states that an infinite simple
group of finite
Morley rank is a linear algebraic group over an algebraically
closed field. In the last fifteen years
an increasing number of partial affirmative answers
have been given to special cases of
this conjecture. Nevertheless
major portions of the problem remain open and
counterexamples are not unexpected.

In recent years, results which elucidate
strong analogies between algebraic groups
and groups of finite Morley rank without
proving specific isomorphism theorems have
reappeared in the area. Such theorems
are reminiscent of the early
work of Daniel Lascar and Bruno
Poizat, later developed by Frank Wagner.
The most important examples are in \cite{frecjalcar},
\cite{jacart}, \cite{bbc} and \cite{bcsstorsion}.

An interesting common point of the main
theorems in these four papers is that they
permit the introduction of an abstract
notion of {\em Weyl group} which corresponds
in the algebraic category to the usual
Weyl group. Several definitions have been suggested
and especially \cite{jacart}
and \cite{bcsstorsion} provide
keys to the relationships among these possibilites.
In this context the centralizers of {\em tori}
are of importance. The following is known for
all connected linear algebraic groups:

\begin{fact}[{\cite[Theorem 22.3]{hu}}]
Let $G$ be a connected linear algebraic group and
$S$ be a torus in $G$. Then $C_G(S)$ is connected.
\end{fact}

The thrust of the present work is an analogous
result for the centralizers
of {\em decent tori}
in connected groups of finite Morley rank:

\begin{theorem}
\label{cdtcon}
Let $T$ be a decent torus of a connected group $G$ of finite
Morley rank. Then $C_G(T)$ is connected.
\end{theorem}

Theorem \ref{cdtcon} will then be used to obtain several
results of varying degrees of generality each of which,
in its own way,
is in the spirit of the aforementioned recent
developments.
Our final and more decisive conclusions, denoted theorems,
in each of the three final sections of this
work depend
on minimality conditions which will be explained
more in detail in Section \ref{minbackg}.

In some cases, such minimality conditions are natural
because $PSL_2$ is the only simple algebraic group where
the result holds, as occurs in Section \ref{weylsection}'s
results on small Weyl groups.
There are other cases however where one expects
more general results to be true and provable.
Our second main result belongs to this latter category:

\begin{theorem}
\label{borelcovering}
Let $G$ be a connected non-solvable group of finite Morley
such that if $A$ is a non-trivial definable
solvable subgroup of $G$
then $N_G^\circ(A)$ is solvable.
Then every element of $G$
belongs to a Borel subgroup. 
\end{theorem}
The analogous
result is known for all connected linear algebraic groups over
algebraically closed fields:

\begin{fact}
[{\cite[Theorem 22.2]{hu}}]
Let $G$ be a connected linear algebraic
group and $B$ a Borel subgroup of $G$.
Then the union of all conjugates of $B$
is $G$.
\end{fact}

In Section \ref{weylsection}, we initiate an analysis
of Weyl groups in relatively general settings and then
draw conclusions on minimal groups. There, minimality
is a natural condition under which the
structure of the Weyl group is very restricted
although still far from being
exactly understood. Theorem \ref{cyclic_cor}
is the final conclusion of Section
\ref{weylsection}:

\begin{theorem}\label{cyclic_cor}
Let $G$ be a minimal connected simple
group of finite Morley rank,
and let $T$ be a decent torus of $G$.
Then the Weyl group
$W= N(T)/C(T)$ is metacyclic.
\end{theorem}
The proof of this theorem involves
steps of varying degrees of generality which
seem important to us for their own sake.
Hence, as the reader of Section \ref{weylsection}
will note, we have taken special care to
extract the important steps of the main proof
as separate conclusions.

Not only do such theorems as Theorem \ref{cyclic_cor}
provide a more precise
description of general structural properties
of minimal groups, but also they are expected to
be of use in classification results in the spirit
of the Cherlin-Zil'ber conjecture. This is an
occurrence of the phenomena in the analysis of
groups of finite Morley rank that have been
witnessed more frequently since \cite{bbc},
and in that sense the present work is a contribution
to the aforementioned line initiated by the work
of Poizat and Lascar.

Our final section is devoted to definable groups
of automorphisms of groups of finite Morley rank.
We obtain general structural information about
automorphisms of minimal connected simple groups.
For reasons which we detail at the beginning of that section,
we feel the conclusions in Section \ref{autosection}
do not only provide specific structural descriptions but also
suggest an appraoch to the Cherlin-Zil'ber conjecture in odd type.
We prove the following theorem:

\begin{theorem}\label{autthm_torii}
Let $\hat{G} = G A$ be a group of finite Morley rank where
G is a definable minimal connected simple subgroup
normal in $\hat G$,
$A$ is a connected definable abelian subgroup, and
$\hat{G}$ is centerless.
The following two properties are contradictory.
\begin{itemize}
\item  $G$ has a divisible abelian torsion subgroup $Q$
 with a non-trivial Weyl group $W= N_G(Q)/C_G(Q)$.
\item  $\hat{G} / G$ has nontrivial divisible torsion.
\end{itemize}
\end{theorem}
On the way to this final theorem, we will prove a more technical
statement which suggests an inductive approach to the same
result without the minimlity hypothesis.

\section{Centralizers of tori}\label{mainsec}

In this section, we prove Theorem \ref{cdtcon}.
The proof necessitates relatively little background
beyond basic stable group theory
or from
the more specialized theory of groups of finite
Morley rank. We will use some basic results
and notions (types, compactness, eq)
from model theory without specific
reference. Poizat's books
\cite{cours} and \cite{grst} are excellent
sources for model theory in general and the
bases of the theory of stable groups. Wagner's
book \cite{stgr} provides technical improvements.
As for groups of finite Morley rank,
Alexandre Borovik and Ali
Nesin's book \cite{bnbook}
is our general reference.

The central concept of this section as well as
most of the sequel is the notion of a {\em decent
torus}. To define it, we need the notion
of the {\em definable hull} of a set in a group
of finite Morley rank. For an arbitrary subset
$X$ of a group $G$ of finite Morley rank,
the definable hull of $X$ with respect to
$G$, denoted $d(X)$,
is the smallest definable subgroup of $G$
which contains $X$. Its existence is a consequence
of the {\em descending chain condition} satisfied
by the definable subgroups of a group of finite Morley
rank. The definable hull of an arbitrary subgroup
of a group of finite Morley provides an extension
of the notion of connected component. Indeed,
if $X$ is an arbitrary subgroup of a group
$G$ of finite Morley rank then $X^\circ$ is
defined as $X\cap d^\circ(X)$.

A {\em decent torus} in a group of finite Morley
rank is a definable, divisible, abelian
subgroup which is the definable hull of its
torsion. This notion coincides with
the usual notion of torus in algebraic
groups over algebraically closed fields
without additional definable subgroups.
In nonzero characteristic,
this behavior of multiplicative subgroups
of algebraically closed fields holds true
in the more general context of fields of finite
Morley rank as well:

\begin{fact}[\cite{wagbadp}]\label{dtinfieldp}
Any definable subgroup of the multiplicative
group of a field of finite Morley rank
of nonzero characteristic is a decent torus.
\end{fact}

More generally,
any group of finite Morley rank which contains
a nontrivial divisible abelian $p$-subgroup
(a {\em $p$-torus})
for an arbitrary
prime $p$ contains a nontrivial
decent torus. Indeed,
it suffices to consider the definable hull of
such a $p$-torus. This strong tie
between divisible abelian torsion
groups and their definable hulls in
groups of finite Morley rank allows
one to transfer various results
on divisible abelian torsion groups
to decent tori. We remind the one
about the {\em rigidity} of tori which will
be used in the sequel without reference:

\begin{fact}[{\cite{BP}}]
Let $G$ be a group of finite Morley rank
and $T$ be a decent torus in $G$.
Then $[N_G(T):C_G(T)]<\infty$ and this
index is bounded.
\end{fact}
This result was proven for $p$-tori
and it is immediate for decent tori
thanks to the following elementary
property of definable hulls in a group
$G$ of finite Morley rank:
$C(X)=C_G(d(X))$ for any $X\subseteq G$.

The following theorem of Cherlin
is fundamental:

\begin{fact}[{\cite{goodtori}}]
\label{dtgen}
Let $G$ be a group of finite Morley
rank. 

1. If $T$ is a decent torus of $G$ then
$\bigcup_{g\in G}\ C^\circ(T^g)$ is generic
in $G$.

2. The maximal decent tori in $G$ are conjugate.
\end{fact}

A definable subset $X$ of a group of finite Morley
rank $G$ is said to be {\em generic} in $G$
if finitely many translates of $X$ cover
$G$, equivalently if $\rk(X)=\rk(G)$, where
$\rk$ denotes the Morley rank.
A definable subset $X$ of a group of finite Morley
rank $G$ is said to be {\em generous}
if $\bigcup_{g\in G} X^g$ is generic
in $G$.
The following fact is a special case of
Lemma 6.2 in \cite{grst}.

\begin{fact}[{\cite[Lemma 6.2]{grst}; see
also {\cite[Lemma 1.4]{bcsstorsion}}}]
\label{genericinquotient}
Let $G$ be a group of finite Morley rank,
$A\subseteq G$ a set of parameters,
$H$ an $A$-definable normal subgroup of $G$
and $\overline G=G/H$. If $a\in G$
is generic over $A$ then its image
$\overline a$ in $\overline G$
is generic over $A$.
\end{fact}
By an element generic over a set in
a group of finite Morley rank $G$, we mean
a realization of a {\em generic type}
of $G$ over this set of parameters. A type
is generic if and only if it contains
only generic sets over the set of parameters
of the type.

Before we start the proof of Theorem \ref{cdtcon}
we also quote two pieces of torsion folklore which
have considerable consequences on the behaviour
of torsion elements in groups of finite Morley rank.

\begin{fact}[{\cite[Ex.~11 p.~93; Ex.~13c p.~72]{bnbook}}]\label{nodeftorsion}
Let $G$ be a group of finite Morley rank and let $H \normal G$ be a
definable subgroup.  If $x\in G$ is an element such that $\bar{x}\in G/H$
is a $p$-element for some
prime $p$, then $x H$ contains a $p$-element.  Furthermore,
if $H$ and $G/H$ are $p^\perp$-groups, then $G$ is a $p^\perp$-group.
\end{fact}

\begin{fact}[{\cite[Ex.~10 p.~93]{bnbook}}]
\label{defnblcyclicstruct}
Let $G$ be a group of finite Morley rank
and $g\in G$. Then $d(g)=A\oplus C$ where
$C$ is finite cyclic and $A$ is divisible.
\end{fact}

We first prove two lemmas of independent interest. The
first one is very much in the spirit of Lemma 1.5
of \cite{bcsstorsion}.

\begin{lemma}\label{da_contains_CTa}
Let $G$ be a group of finite Morley rank and
 take $a \in G$ generic in some coset of $G^\circ$.
Suppose that $G$ normalizes a maximal decent
torus of $C(a)$.
Then $d(a)$ contains this decent torus.
\end{lemma}
\proof
We may assume $G=d(a)G^\circ$. Let $T$ be the
maximal normal decent torus in $G$.
Let $T_0$ be the torsion in $T$ and let $R = d(a) \cap T$.
As $G^\circ$ centralizes $T$,
 $R$ does not depend upon $a$, but merely upon $a G^\circ$.
Thus the definition of $R$ does not necessitate
additional parameters to those needed to define $a G^\circ$.
So, by Fact \ref{genericinquotient},
 the image $\bar a$ of $a$ is generic in $\bar G = G/R$.
We may therefore assume that $d(a) \cap T = R = 1$.
So it suffices to show that $S_0 = C_{T_0}(a)$ is finite.

By Fact \ref{defnblcyclicstruct}, $d(a)=A \oplus C$ is the direct
sum of a divisible abelian group $A$ and a finite cyclic group $C$.
If $n=|C|$ then for any multiple $N$ of $n$, we have $d(a^N)=A$.
On the other hand, for any torsion element $s \in S_0$,
since $T_0\subseteq\acl(\emptyset)$,
 the element $a'=as$ is also generic over $\emptyset$.
Thus, being in the same coset of $G^\circ$,
$a'$ and $a$ realize the same type.
Letting $N$ be a multiple of $n$ and the order of $s$,
 it follows that $d((a')^n)=d((a')^N)=d(a^N)=d(a^n)$
 and thus $s^n \in d(a^n) \le d(a)$.
Now by our reductions $d(a)$ contains no $p$-torus for any $p$,
 and hence $d(a)$ has bounded exponent.
Thus $s^n$ has bounded exponent, with $s$ varying and $n$ fixed,
 and so $S_0$ is finite as claimed.
\qed

The second preparatory lemma
is reminiscent of Fact \ref{genericovering} below.

\begin{lemma}\label{genericoveringtype}
Let $G$ be a group of finite Morley rank.
Let $H$ be a definable subgroup of $G$ and
$X$ a definable subset of $G$ such that the
Morley degree of $X$ is $1$, and
\[\rk(X\cap (\bigcup_{g\notin H}X^g))=\rk(X).\]
If $x$ is a realization of the generic
type of $X$ then there exists
$g\in G\setminus H$ such that $x^g$ also realizes
the generic type of $X$.
\end{lemma}
\proof
Let $\Phi(x)$ be the generic type of $X$
over the parameters for $H$ and $X$.
We introduce the following set of
formulas on the parameters
for $H$ and $X$:
\[\Psi(x,y)=\{y\not\in H\}\cup
\{X(x), X(x^{y^{-1}}): X(x)\in\Phi(x)\}\]
Here $X(x)$ ranges over the entire
set of formulas in $\Phi(x)$.

We claim that
$\Psi(x,y)$ is consistent.
Any finite subset of $\Psi(x,y)$ can
be assumed to be of the following kind
\[X_1(x)\wedge X_2(x^{y^{-1}})\wedge
y\not\in H\]
where $X_1(x)$ and $X_2(x)$ are in
$\Phi(x)$. Moreover, since $\Phi(x)$
is the generic type of $X$ which is of Morley
degree $1$, replacing $X_1(x)$ with
$X_1(x)\wedge X_2(x)$ and
$X_2(x^{y^{-1}})$ with
$X_1(x^{y^{-1}})\wedge X_2(x^{y^{-1}})$
we may assume that we are dealing
with $X(x)\wedge X(x^{y^{-1}})$
for some $X(x)\in\Phi(x)$.
By the generic covering assumption,
there exist $b\not\in H$ and
$a,a_1\in X$ such that $a_1^b=a$.
Hence, $(a,b)$ satisfies the formula
in question.

By compactness, $\Psi(x,y)$
is realized in an elementary extension
of $G$ by $(\alpha,\beta)$. This means
that both $\alpha$ and $\alpha^{\beta^{-1}}$
realize $\Phi(x)$.
It follows
that $\alpha$ is generic in $X$
and $X^\beta$.
\qed

 In the proof of Theorem \ref{cdtcon},
we will need several facts on genericity.

\begin{fact}[{\cite[Lemma 2.2]{bcsstorsion}}]
\label{genericgeneration}
Let $G$ be a group of finite Morley rank, $H$
a connected definable subgroup, and $X$ a definable
generic subset of the coset $aH$. Then
$\langle X\rangle=\langle aH\rangle=
\langle a,H\rangle$.
\end{fact}

\begin{fact}[{\cite[Lemma 4.1]{bbc}}]
\label{genericovering}
Let $G$ be a group of finite Morley rank, $H$
a definable subgroup of $G$, and $X\subseteq G$
a definable set such that
\[ \rk(X\setminus\bigcup_{g\not\in H}X^g)\geq
\rk(H)\]
Then $\rk(\bigcup X^G)=\rk(G)$.
\end{fact}

\begin{repstat}{Theorem \ref{cdtcon}}
Let $T$ be a decent torus of a connected group $G$ of finite
Morley rank. Then $C_G(T)$ is connected.\\
\end{repstat}
\proof
We may assume that $T$ is a maximal decent torus
 such that $C(T) > C^\circ(T)$.
Let $H = C^\circ(T)$.
There is a $p$-element $w \in C(T) \setminus H$
 by Fact \ref{nodeftorsion}.

Consider a realization $a$ of the generic type $\phi(x)$ of $w H$
 over the parameter set $T\cup \{w\}$.
As $w,H \leq C(T)$, clearly $a \leq C(T)$.
For any decent torus $R$ of $C(a)$ containing $T$,
 we have $a \notin H \geq C^\circ(R)$, and $R = T$ by maximality of $T$.
So $T$ is a maximal decent torus of $C(a)$.
As $T$ is normal in $H$,
 $T \leq d(a)$ by Lemma \ref{da_contains_CTa}.

Before moving ahead, we make a remark about
the setwise stabilizer of $wH$, namely the definable
subgroup of $G$:
\[\Stab(wH)=\{g\in G:\rk(wH\triangle (wH)^g)=\rk(wH)\}.\]
By Fact \ref{genericgeneration}, $\Stab(wH)\leq
N_G(\langle w, H\rangle)\leq N_G(H)$, the last
inclusion being justified by the fact that
$H=\langle w,H\rangle^\circ$. On the other hand,
$N_G^\circ(H)=C^\circ(T)=H\leq\Stab(wH)$. Thus,
$\Stab^\circ(wH)=H$.

Suppose first that the set
$wH\setminus\bigcup_{g\notin\Stab(wH)}(w H)^g$
is generic in $wH$.
Then, by the remark in the preceding paragraph,
$\rk(wH\setminus\bigcup_{g\notin\Stab(wH)}(w H)^g)
=\rk(H)=\rk(\Stab(wH))$. By Facts \ref{genericovering}
and \ref{dtgen} (1),
a realization $a$ of $\phi(x)$ lies inside
 $C^\circ(S)$ for some maximal decent torus $S$.
As $T \leq d(a) \leq C^\circ(S)$, we have $T \leq S$.
So $a \in C^\circ(S) \leq H$, contradicting $w \notin H$.
Hence $wH\setminus\bigcup_{g\notin\Stab(wH)}(w H)^g$
is not generic in $wH$.

The preceding conclusion implies that
$wH\cap(\bigcup_{g\notin\Stab(wH)}(wH)^g)$
is generic in $wH$. 
For $g\in N_G(H)\setminus\Stab(wH)$ then
$(wH)^g\cap wH=\emptyset$.
It follows
that the set $wH\cap(\bigcup_{g\notin N_G(H)}(wH)^g)$
is generic in $wH$ as well.

At this point, we apply Lemma \ref{genericoveringtype}
with $X=wH$ and $H$ (in the statement of lemma)
for $N_G(H)$. This yields $a$ a generic of $wH$,
and $b\notin N_G(H)$ such that $a^b$ is also generic
in $wH$.
Thus, $T \leq d(a)$ and
$T^{b^{-1}} \leq d(a)$ by Lemma \ref{da_contains_CTa}.
As $d(a)$ is abelian and $T$ is maximal
in $C(a)$, $T^b = T$, and $b \in N(T)\leq N_G(H)$,
 a contradiction to the choice of $b$.
\qed

\section{Background on minimal groups}
\label{minbackg}

As was mentioned in the introduction, many
results proven in the rest of the paper, in particular
the final conclusions of the main sections,
depend on minimality conditions. The following
hypothesis describes the most general minimal
setting in which we will be working in the sequel.

\begin{hypothesis}\label{ngroup}
$G$ is a connected non-solvable group of finite Morley
such that if $A$ is a non-trivial definable
solvable subgroup of $G$
then $N_G^\circ(A)$ is solvable.
\end{hypothesis}
In literature, one generally encounters
a stronger minimality assumption which is
the {\em minimal simplicity}, namely that $G$
is a connected simple group of finite Morley rank
of which the proper definable connected subgroups
are solvable. In Section \ref{weylsection},
we will restrict our analysis to connected minimal
simple groups, while in Section \ref{autosection}
the setting described by Hypothesis \ref{ngroup}
will provide the right context to analyze the
automorphisms of minimal connected simple groups.

It is worth noting that Adrien Deloro and Eric
Jaligot have undertaken a project which investigates
possible generalizations of all results known
under Hypothesis \ref{ngroup} to a wider setting
where the underlying assumption replaces
$A$ with $A^\circ$. They refer to groups
satisfying the weaker hypothesis as {\em
locally$^\circ$ solvable$^\circ$} while the
stronger hypothesis referred to by dropping
the appropriate $\circ$ sign.
We will not use their terminology
here to avoid confusion with
the traditional use of the phrase ``locally solvable''.

The main objective of this section is to review
some of the results on minimal groups
needed in the sequel. Several of these results
were proven in the context of minimal
simple groups. For the sake
of completeness, we will provide their proofs
using Hypothesis \ref{ngroup}
although there are practically
no changes. This will be an opportunity to revisit
them and occasionally
extract additional information from
these proofs which is somewhat hidden in their
form in print.

The following lemma about the intersections
of Borel subgroups is of crucial importance.
A {\em Borel subgroup} of a group of finite
Morley rank is a maximal, definable, connected
solvable subgroup.
Lemma \ref{unipjal} is adapted from \cite{burbender} where
it was proven for minimal connected simple groups.

\begin{lemma}[{\cite[Lemma 2.1]{burbender}}]
\label{unipjal}
Let $G$ be a group of finite Morley rank
satisfying Hypothesis \ref{ngroup}. Let $B_1$ and
$B_2$ be two distinct Borel subgroups satisfying
$U_{p_i}(B_i)\neq1$ where the $p_i$ are
two prime numbers $(i=1,2)$ which are not necessarily
distinct.
Then $F(B_1)\cap F(B_2)=1$.
\end{lemma}
For any group $H$ of finite Morley rank and
prime number $p$, $U_p(H)$
denotes its largest, definable, connected, normal,
nilpotent $p$-subgroup. In general, definable,
connected, nilpotent subgroups of bounded
exponent of groups of finite Morley rank are
called {\em unipotent} subgroups.

We recall two important facts about unipotent
$p$-subgroups of groups of finite Morley rank. They will
be used in the proof of Lemma \ref{unipjal}.

\begin{fact}[{\cite{N1}}]\label{Upnilpotence}
Let $H$ be a connected solvable group of finite Morley rank.
Then there is a unique maximal unipotent $p$-subgroup $U_p(H)$ of $H$,
and $U_p(H) \leq F^\circ(H)$.
\end{fact}

\begin{fact}[{\cite{N3}}]
\label{Ucenter}
Let $G$ be a nilpotent group of finite Morley rank satisfying
$U_p(G) \neq 1$. Then $U_p(Z(G)) \neq 1$.
\end{fact}

\proofof{Lemma \ref{unipjal}}
We first show that $U_p(B_1 \cap B_2) = 1$ for an arbitrary
prime number $p$.
Let  $X = U_p(B_1 \cap B_2)$ and suppose toward
a contradiction that $X \neq 1$.
We may assume that $\rk(X)$ is
maximal among all choices of $B_1$ and
$B_2$.  By Hypothesis \ref{ngroup},
there is a Borel subgroup $B$ of $G$
containing $N^\circ_G(X)$.

We will arrive at a contradiction
by showing that $B_i = B$ for both $i=1,2$.
If $X<U_p(B_1)$ then Fact \ref{Upnilpotence}
implies that
\[X < U_p(N_{U_p(B_1)}^\circ(X)) \leq U_p(B\cap B_1). \]
It follows from the maximal choice of $X$ that $B=B_1$.
Applying the same argument to $B_2$ implies
that $X=U_p(B_2)$ since it has been assumed
that $B_1\neq B_2$. But then,
$B_2\leq N^\circ(X)\leq B\cap B_2$ and thus
$B=B_2$, a contradiction. Thus,
$X=U_p(B_1)$. A symmetric argument
shows that $X=U_p(B_2)$. As a result
$B_1=B=B_2$, and this contradiction
forces $X=1$.

Having proven that $U_p(B_1 \cap B_2) = 1$
for any prime $p$,
we proceed to prove the lemma.
Suppose toward a contradiction that there
is an $f \in F(B_1) \cap F(B_2)$ with $f \neq 1$.
Let $Z_i=U_{p_i}(Z(F(B_i)))$.
By  Fact \ref{Ucenter}, $Z_i$
is a nontrivial unipotent $p_i$-subgroup
of $Z(F(B_i))$ for $i=1,2$.
Let $B$ be a Borel subgroup containing $C^\circ_G(f)$.
As $Z_i \leq Z(F(B_i)) \leq C^\circ_G(f)$
for $i=1,2$, we find
$Z_i \leq U_{p_i}(B_i \cap B)$ for both
$i=1,2$. Thus $B_1 = B = B_2$ by the first part.
\qed

As a next step in our preparation
we adapt a result from \cite{bbc}
to the context of groups
satisfying Hypothesis \ref{ngroup}.
This is an occasion for rewording
the statement of this fact and
extracting some useful information
from its proof.

\begin{fact}[{\cite[Proposition 8.1]{bbc}}]
\label{weakesgenmin}
Let $G$ be a minimal, connected, simple group of
finite Morley rank. Then the set of elements of $G$
which belong to some definable
connected nilpotent subgroup of
$G$ contains a definable generic subset of $G$.
\end{fact}

The proof of Fact \ref{weakesgenmin} yields
in fact the following proposition:
\begin{prop}
\label{BBC_8.1}
Let $G$ be a group of finite Morley
rank satifying Hypothesis \ref{ngroup}.
Then
\begin{itemize}
\item[-] either $G$ does not have torsion,
\item[-]or $G$ has a generous Carter subgroup.
\end{itemize}
\end{prop}

Before proving Proposition \ref{BBC_8.1},
we go over some important notions and facts
needed for its proof as well as for later use.
A {\em Carter subgroup} of a group of finite
Morley rank is a definable, connected, nilpotent
subgroup of finite index in its normalizer.

\begin{fact}\label{classicarter}\

1. {\em (\cite{frecjalcar})} Every group
of finite Morley rank has a Carter
subgroup. The centralizer of a maximal
decent torus contains a Carter subgroup
of the ambient group.

2. {\em (\cite{carwag})} Carter subgroups
of solvable connected groups of finite
Morley rank are conjugate.

3. {\em (\cite{frec1})} Carter subgroups
of connected solvable groups of finite
Morley rank are connected.

4. {\em (\cite{frec1})} Carter subgroups
of connected solvable groups of finite Morley
rank are self-normalizing.

5. {\em (\cite{frecjalnewt})} Carter subgroups
of connected solvable groups of finite Morley
rank are generous.
\end{fact}

For genericity arguments, we will need the following
fact in addition to
what was recalled in Section \ref{mainsec}:

\begin{fact}
\label{jalsnfolklore}

1. {\em {\cite[Lemma 2.2]{jacart}}}
A generous subgroup of a group of finite Morley
rank is of finite index in its normalizer.

2. {\em \cite[Lemma 2.4]{jacart}}
Let $G$ be a group of finite Morley rank
and $H$ a connected
generous subgroup of $G$.
Then any definable generic subset $X$ of $H$
is generous in $G$.
\end{fact}

The analysis of torsion elements in the proof
of Proposition \ref{BBC_8.1} uses the
following recent theorem:

\begin{fact}[{\cite[Theorem 4]{bbc},\cite{bcsstorsion}}]
\label{notpdeg}
Let $G$ be a connected group of finite
Morley rank, and $a\in G$ a $\pi$-element
where $\pi$ is a set of primes.
Then $C^\circ(a)$
contains an infinite abelian $p$-subgroup
for some $p\in\pi$.
\end{fact}

The following lemma is extracted from the proof
of Fact \ref{weakesgenmin} in \cite{bbc}:
\begin{lemma}\label{unipotence_genericity}
Let $G$ be a group of finite Morley rank
satisfying Hypothesis \ref{ngroup}.
Suppose that $G$ contains no divisible torsion
but has nontrivial unipotent torsion.
Let $B$ be a Borel subgroup of $G$ with $U_p(B) \neq 1$.
Then $\bigcup_{g\in G} B^g$ is generic in $G$.
\end{lemma}
\proof
By assumption,
$G$ has no nontrivial decent tori but has nontrivial
unipotent $p$-subgroups. Then Hypothesis \ref{ngroup}
and Fact \ref{Upnilpotence} imply that
$G$ has a Borel subgroup $B$ with $U_p(B)\neq1$.
Note that Fact \ref{nodeftorsion}
and the fact that $G$ does not have decent tori
imply that $B/U_p(B)$ does not have $p$-torsion. Using
Fact \ref{dtinfieldp} and Zil'ber's field interpretation
results in conjunction with Fact \ref{nodeftorsion},
one deduces that $B=U_p(B)C_B(U_p(B))$.

Next, one proves that for a generic element $g\in B$,
$d(g)\cap U_p(B)\neq1$.
Suppose otherwise. Let $z$ be an element
of order $p$ in $Z(U_p(B))$.
Then $zg$ is still a generic in $G$.
We then have $(zg)^p=z^pg^p=g^p\in
d(zg)\cap d(g)$. By the assumption on $d(g)$ and
Fact \ref{nodeftorsion}, $d(g)\leq d(zg)$ and it
follows that $z\in d(g)$. This contradicts
the assumption on $d(g)$.

It follows that $B$ is generically disjoint from its conjugates
as otherwise we would have distinct Borel subgroups meeting
$U_p(B)$ nontrivially, contradicting Lemma \ref{unipjal}.
Since $[N_G(B):B]<\infty$, we can
apply Fact \ref{genericovering} with $H=N_G(B)$ and $X=B$.
Thus, $B$ is generous in $G$.
\qed

Now, we proceed to the proof of Proposition \ref{BBC_8.1}.\\

\proofof{Proposition \ref{BBC_8.1}}
We may assume that $G$ has $p$-torsion for some
fixed prime $p$. Then by Fact \ref{notpdeg},
$G$ has infinite abelian $p$-subgroups.

We first analyze the case when $G$ has nontrivial
decent tori. Let $T$ be a nontrivial decent
torus in $G$. By Fact \ref{dtgen}, $C^\circ(T)$
is generous in $G$. By Hypothesis \ref{ngroup},
$C^\circ(T)$ is contained in a Borel subgroup $B$ of
$G$. In particular, $B$ is generous in $G$.
By Fact \ref{classicarter}.1, $B$ has Carter
subgroups, and
by Fact \ref{classicarter}.5, a Carter subgroup $C$ of $B$
is generous in $B$. By Fact \ref{jalsnfolklore}.2, $C$ is generous
in $G$ and it follows from Fact \ref{jalsnfolklore}.1,
that $[N_G(C):C]<\infty$.

We may assume by the previous paragraph that
$G$ has no nontrivial decent tori. As a result
of this and Fact \ref{notpdeg}, $G$ has nontrivial
unipotent $p$-subgroups. Hence, we can apply
Lemma \ref{unipotence_genericity} to conclude
the presence of a
generous Borel subgroup $B$.
As in the case when $G$ contained nontrivial
decent tori, the Carter subgroups of $B$
are generous Carter subgroups of $G$.
\qed

\section{Covering by Borel subgroups}

In this section, we prove Theorem \ref{borelcovering}.
As was explained in the introduction, one expects
Theorem \ref{borelcovering} to hold for all
connected groups of finite Morley rank.
Here we are content with those satisfying
Hypothesis \ref{ngroup}. The two technical
steps of the proof (Lemma \ref{unipotentselfnormalization}
and
Corollary \ref{pcentbor}) seem to be of independent
interest.

The following recent theorem is essential for the
sequel:

\begin{fact}
[{\cite[Theorem 3, Corollary 3.3]{bcsstorsion}}]
\label{torictorsion}
Let $G$ be a connected group of finite
Morley rank, $\pi$ be a set of primes, $a$ any
$\pi$-element of $G$, and suppose
$C_G^\circ(a)$ has $\pi^\perp$-type.
Then $a$ belongs to any maximal $\pi$-torus
of $C(a)$.
\end{fact}
Given a set $\pi$ of prime numbers,
a group of finite Morley rank is said to
be of {\em $\pi^\perp$-type} if it does not
contain nontrivial unipotent $p$-subgroups
for any $p\in\pi$.

We also need a well-known fact about
the structure of Hall $\pi$-subgroups
of connected solvable groups of finite
Morley rank:

\begin{fact}[{\cite[Theorem 9.29]{bnbook};
\cite[Corollaire 7.15]{frec1}}]
\label{consolvhall}
Hall $\pi$-subgroups of a connected solvable
group of finite Morley rank are connected.
\end{fact}

One of the main ingredients of the proof
is the following lemma which seems valuable in
its own right.

\begin{lemma}\label{unipotentselfnormalization}
Let $G$ be a group
of finite Morley rank satisfying
Hypothesis \ref{ngroup}
and $B$ be a Borel subgroup
of $G$ such that $U_p(B)\neq1$ for some prime
number $p$. Then $p\not|[N_G(B):B]$.
\end{lemma}

\proofof{Lemma \ref{unipotentselfnormalization}}
 Suppose towards a contradiction that
$p\ |[N_G(B):B]$ and that this is witnessed
by an element $x$ in $N_G(B)\setminus B$.
Using Fact \ref{nodeftorsion},
we may assume that $x$ is a $p$-element
whose order modulo $B$ is $p$. Let $U$ denote
the maximal unipotent $p$-subgroup of $B$.
By assumption $U\neq1$.
We will reach a contradiction
using the generic behavior of $xB$.

We will apply Fact \ref{genericovering}
with $X=xB$ and $H=N_G(B)$.
Let $X_1= xB\cap\left(\bigcup_{g\not\in N_G(B)}(xB)^g
\right)$.
If the set $X_1$ is generic
in $xB$ then there exists
$y\in xB\cap (xB)^g$
for some $g\not\in N_G(B)$. In particular,
$y^p\in B\cap B^g$.
Since $x$ normalizes
$B$ and $x^g$ normalizes
$B^g$, $y$ normalizes both $B$ and $B^g$. Considering
$d(y)\cap B\cap B^g$, one can translate $y$
to a $p$-element $z$ by an element of $B\cap B^g$
(Fact \ref{nodeftorsion}).
It follows that $z$ centralizes infinite
subgroups of both $U$ and $U^g$. This
forces $g\in N_G(B)$ again by 
Lemma \ref{unipjal}. This is a contradiction.

If the set $X_1$ is not generic in $xB$
then by Fact \ref{genericovering}
the set $\bigcup (xB)^G$ is generic in $G$.
At this point we have to consider two possible
cases depending on whether or not $G$ has a
nontrivial decent torus.

First, assume that $G$ has a nontrivial decent
torus $T$ which we may assume to be maximal.
By Fact \ref{dtgen}, there exist $y\in xB$ and $g\in G$
such that $y\in C^\circ(T^g)$ by
Fact \ref{dtgen}. Now we can translate
$y$ by an element of $d(y)\cap B$ to a $p$-element
$y'\in xB$ (Fact \ref{nodeftorsion}). Since $y'\in N_G(B)$,
$C^\circ(y')$ contains a nontrivial unipotent
$p$-subgroup of $B$. Hence, $C^\circ(y')$
is contained in $B$ by Lemma \ref{unipjal}.

Since the
translation to $y'$ is done by an element
of $d(y)$, $y'$ stays in
$C^\circ(T^g)$. Let $B_2$ be a Borel
subgroup of $G$ such that
$y'\in C^\circ(T^g)\leq B_2$
(Hypothesis \ref{ngroup}). Since $y'$
is a $p$-element, it is contained in a
Sylow $p$-subgroup $P$ of $B_2$. Since
$B_2$ is connected, Fact \ref{consolvhall}
implies that $P$ is connected.
If $P$ is a $p$-torus, then
$y'\in P\leq C^\circ(y')\leq B$. This
contradicts that $y'\not\in B$. Thus
$B_2$ contains nontrivial
unipotent $p$-torsion. Hence, so does
$C_{B_2}^\circ(y')$. It follows from
Lemma \ref{unipjal} that
$B_2=B$, which again forces
$y'\in B$, a contradiction.

Second, we assume that $G$ has no
divisible torsion. Nevertheless,
$G$ has $p$-torsion. Thus 
by Fact \ref{notpdeg}, $G$ has unipotent
$p$-torsion. By Lemma \ref{unipotence_genericity},
a Borel subgroup which contains
a nontrivial unipotent $p$-subgroup is generous.
In our situation $B$ is one such. Hence
there exist $y\in xB$ and $g\in G$ such that
$y\in B^g$. Note that $g\not\in N_G(B)$
since $y\not\in B$. Since $y\in B^g$,
so is $y^p$. Since $y^p\in B$ as well,
we can repeat the argument which
yields a contradiction in the case
where $xB$ is generically covered by
$\bigcup_{g\not\in N_G(B)}(xB)^g$.
No other case remains.
\qed

The following corollary is potentially useful:

\begin{cor}\label{pcentbor}
Let $G$ be a group of finite Morley
rank which satisfies Hypothesis \ref{ngroup}.
If $x$ is an element of finite order
of $G$ then $x$ lies inside
any Borel subgroup $B$ containing $C^\circ(x)$.
\end{cor}
\proof By way of contradiction,
suppose that $x$ is an element of
finite order $m$ for which there
exists a Borel subgroup $B$ of $G$ such that
$x\not\in B$ but $C^\circ(x)\leq B$.
Choose $x$ so that $m$ is minimal among
the orders of all such offending elements.
Suppose
$m=p_1^{\alpha_1}\ldots p_k^{\alpha_k}$
with all the $\alpha_i$ different from zero.
We set $\pi=\{p_1,\ldots,p_k\}$,
and for
each $i\in\{1,\ldots,k\}$, we let
$m_i=\frac{m}{p_i^{\alpha_i}}$ and
$x_i=x^{m_i}$. Note that
$x=x_1\ldots x_k$, and
for each $i\in\{1,\ldots,k\}$,
$|x_i|=p_i^{\alpha_i}$.

If $C_G^\circ(x)$ is of $\pi^\perp$-type
then it follows from
Fact \ref{torictorsion}
that $x\in C_G^\circ(x)$. Thus,
$C_G^\circ(x)$ contains nontrivial unipotent
$p_i$-torsion
for some $p_i\in\pi$.
Then by Lemma \ref{unipjal} and the standing
assumption on $x$, $x\in N_G(B)\setminus B$.
Moreover, since $C_G(x)\leq C_G(x_i)$,
it follows from Lemma \ref{unipjal} that
$C_G^\circ(x_i)\leq B$. Then, by
Lemma \ref{unipotentselfnormalization},
$x_i\in B$. Thus, $|\pi|>1$. As a result,
we can apply induction to the other $x_j$.
Since $C_G(x)\leq C_G(x_j)$ for every
$j\in\{1,\ldots,k\}$, it follows as for $x_i$
that $C_G^\circ(x_j)\leq B$.
The minimality of $m$ forces that $x_j\in B$.
Thus, $x\in B$, a contradiction.
\qed
We note that $a$ need not lie inside $C^\circ(a)$
if $G$ is a bad group whose Borel subgroups are nilpotent
groups constructed by Baudisch (\cite{baudgroup}).

Let us restate Theorem \ref{borelcovering}
before giving its proof which is now immediate.\\

\begin{repstat}{Theorem \ref{borelcovering}}
Let $G$ be a group of finite Morley
rank which satisfies Hypothesis \ref{ngroup}.
Then every element of $G$
belongs to a Borel subgroup.\\ 
\end{repstat}
\proof
Let $x$ be an offending element. By
Corollary \ref{pcentbor}, we may assume
$x$ is of infinite order.
Fact \ref{defnblcyclicstruct} yields
an $x_1$ such that
$d(x)=\langle x_1\rangle\oplus
d^\circ(x)$. By Corollary \ref{pcentbor},
any Borel containing $C_G^\circ(x_1)$
contains also $x_1$. Since $d^\circ(x)\leq
C_G^\circ(x_1)$, any Borel containing
$x_1$ contains $x$ as well.
\qed

\section{Weyl groups in minimal groups}
\label{weylsection}

We now turn our attention to the smallest Weyl groups.
In the context of groups of finite Morley rank,
there are several possible ways of defining a Weyl
group. Here we will refer to any subgroup
of $N_G(T)/C_G(T)$ where $G$ is a group of finite
Morley rank and $T$ is a $\pi$-torus (equivalently,
a decent torus) as a Weyl group. 
It is worth reminding that $C_G(T)$ is
a connected group by Theorem \ref{cdtcon}
when $G$ is connected.

Before
moving any further, we emphasize that
in a connected group $G$ of finite Morley rank
if $T$ is a minimal decent torus and
$a$ is a representative of a Weyl group
element defined with respect to $T$
then the torsion subgroup of $C_T(a)$ is finite.
We shall find below that this condition also holds
in minimal connected simple groups of finite Morley
rank.
So we now turn our attention to its consequences.

\begin{prop}\label{metacyclic}
Let $G$ be a connected group of finite Morley rank,
 let $T$ be a $\pi$-torus of $G$, and
 consider some $W \leq N(T) / C(T)$.
Suppose that $C_T(a)$ is finite for all $a \in W^\#$.
Then $W$ is a Frobenius complement.
\end{prop}

A {\em Frobenius complement}
is the stabilizer of a point in a {\em Frobenius group},
which is a transitive permutation group
 such that no non-trivial element fixes more than one point
 and some non-trivial element fixes a point.
Such groups have a quite restrictive structure
especially when they are finite (see Fact \ref{Frobenius}
below).

To prove Proposition \ref{metacyclic},
we will need facts from various algebraic
sources. The first one is
about the Tate module construction of a complex
representation.

\begin{fact}[{\cite[\S3.3]{ber1,berbor01}}]\label{endptor}
Let $T$ be a $p$-torus of Pr\"ufer $p$-rank $n$.
Then $\End(T)$ can be faithfully represented
as the ring $M_{n \times n}(\zz_p)$
of $n \times n$ matrices over the $p$-adic integers $\zz_p$.
\end{fact}

The proof of this fact produces an isomorphism
 $\psi : \End(T) \to \invlim_k \End(\Omega_k(T)) \cong M_{n \times n}(\zz_p)$
from the maps $\psi_k : \End(T) \to \End(\Omega_k(T))$
and the sequence 
$$ \End(\Omega_1(T)) \leftarrow \End(\Omega_2(T))
\leftarrow \End(\Omega_3(T)) \leftarrow \cdots $$
So any $w \in \End(T)$ has an infinite kernel if
and only if
it has a kernel in each $\End(\Omega_k(T))$, and then
zero as an eigenvalue in $M_{n \times n}(\zz_p)$.

Some well-known results about the representations
of finite groups are indispensable for our arguments.
Maschke's
theorem is one of them:

\begin{fact}[{\cite[Corollary 1.6]{FuHa}}]\label{irredrep}
Any complex representation of a finite group
is a direct sum of irreducible representations.
\end{fact}

Another representation theory fact which is crucial
for the proof of Proposition \ref{metacyclic}
is Brauer's theorem on {\em splitting fields}:

\begin{fact}[{\cite[Theorem 5.25]{ba2}}]\label{splitfield}
Let $G$ be a finite group of exponent $m$. Then
the cyclotomic field $\Lambda^{(m)}/\qq$ of the $m$th
roots of unity is a splitting field for $G$.
\end{fact}
For further details we refer the reader to Section 5.12 of
\cite{ba2}.

The proof will involve reduction modulo prime numbers:

\begin{fact}[{\cite[Theorem 4.38]{ba1}}]\label{modpreduction}
Let $f(x)\in\zz[x]$ be a monic polynomial of degree
$n$, $E$ a splitting field of $f(x)$ over $\qq$, $p$
a prime such that $f_p(x)$ (the image of $f(x)$ after
reducing modulo $p$) has distinct roots in its splitting
field $E_p$ over $\zz/p\zz$. Let $D$ be the subring
of $E$ generated by the roots of $f(x)$. Then

$($a$)$ There exist homomorphisms $\psi$ of $D$ into
$E_p$.

$($b$)$ Any such homomorphism gives a bijection of the
set $R$ of roots $f(x)$ in $E$ onto the set $R_p$
of roots $f_p(x)$ in $E_p$.

$($c$)$ If $\psi$ and $\psi'$ are two such homomorphisms
then $\psi'=\psi\sigma$ where
$\sigma\in$Gal$(E/\qq)$.

\end{fact}


Last but not least before proving Proposition
\ref{metacyclic} is an important notation
which will be needed in the sequel
as well. For a prime $p$,
a $p$-torus, i.e. a divisible abelian $p$-group,
is the direct sum of copies of $\zz_{p^\infty}$,
the group of complex $p^n$th roots of unity under multiplication.
When the ambient group is of finite Morley rank,
the number of copies is finite for a given $p$-torus
$T$ and it is denoted by $\pr_p(T)$, the {\em Pr\"ufer
$p$-rank} of $T$.\\

\proofof{Proposition \ref{metacyclic}}
Let $T_p$ be the maximal $p$-torus in $T$,
 choosing $p \in \pi$ such that $T_p \neq 1$.
Since $C_T(a)$ is finite for any $a \in W$, $W \hookrightarrow \End(T_p)$.
By Fact \ref{endptor}, there is a faithful representation of $W$
 as $n \times n$ matrices over the $p$-adic integers $\zz_p$
 where $n= \pr_p(T)$.
By tensoring with $\cc$, viewed as the algebraic closure of
 the fraction field $\qq_p$ of $\zz_p$,
we have a complex representation $V$ of $W$.
By Fact \ref{irredrep}, $V = V_1 \oplus \cdots \oplus V_l$
 is a direct sum of irreducible $W$-modules $V_i$.

Suppose that some $a \in W^\#$ has $C_{V_i}(a) \neq 0$.
Then $1$ is an eigenvalue of $a$. 
So $1 - a \in M_{n \times n}(\zz_p)$ has zero as an eigenvalue.
But now $1 - a \in \End(T_p)$ has zero as an eigenvalue.
So $C_{T_p}(a)$ is infinite, a contradiction.

Let $e$ be the exponent of $W$.
Since $W$ is a finite subgroup of $\GL_n(\cc)$,
 $W$ lies inside $\GL_n(\qq[r])$ where $r$ is a
primitive $e$th root of unity (Fact \ref{splitfield}).
We note that $x^e-1$ has distinct
roots over any field with characteristic
relatively prime to $e$
because $x^e-1$ is relatively
prime to its derivative $e x^{e-1}$.
All such roots over $\cc$ clearly lie in $\zz[r]$.
Let $q$ be a prime coprime to $|W|$.
By Fact \ref{modpreduction} (b),
 the map
$\zz[r] \to(\zz/q\zz)[r]$ induces
a bijection on the roots of $x^e-1$.
For any $w \in W^\#$, the characteristic
polynomial $p_w(x)$ has the form
 $\Pi_{i=1}^n (x - \lambda_i)$ where
each eigenvalue $\lambda_i$ is a root of $x^e-1$.
As no $\lambda_i$ is 1 in $\zz[r] \leq \cc$,
 their images $\lambda_i/(q)$ in $(\zz/q\zz)[r]$
under this bijection are also not 1.
In particular $w \neq 1 \pmod{q}$ and 
we still
have a faithful representation of $W$ on $V = ((\zz/q\zz)[r])^n$
 such that again $C_V(a) = 0$ for every $a \in W^\#$.

Now let $F=V\rtimes W$. Then the group $F$
is a transitive permutation group
 in its affine action on $V$.
We show that $F$ is a Frobenius group with $W$
a Frobenius complement.
If some $(b,a) \in F^\#$ fixes two points $x,y \in V$
 then $a x + b = x$ and $a y + b = y$, and $(a-1) (x - y) = 0$.
If $x \neq y$ then 1 is an eigenvalue of $a$,
 so $a = 1$ and $b \neq 0$.
But when $b \neq 0$, it does not fix any element of $V$.
Since
$W \cong (0,W) \leq F$ clearly stabilizes $0 \in V$,
$F$ is a Frobenius group.
For $f \in F \setminus W$,
 $\langle f,W\rangle$ meets $V \cong (V,1) \leq F$.
So $W$ is the whole stabilizer of $0 \in V$,
 and $W$ is a Frobenius complement.
\qed

To apply Proposition \ref{metacyclic}, we need the following
basic facts about the structure
of a Frobenius complement in finite groups.

\begin{fact}[{\cite[10.3.1 p.~339]{gor}}
]\label{Frobenius}
Let $W$ be a Frobenius complement in a finite
Frobenius group.  Then
\begin{enumerate}
\item Any subgroup of $W$ of order
$pq$, $p$ and $q$ primes, is cyclic.
\item Sylow $p$-subgroups of $W$ are either cyclic or
 possibly generalized quaternion if $p=2$.
\item If $W$ is even, $W$ possesses a unique
involution which is central.
\end{enumerate}
\end{fact}

As it is clearly visible in the statement
of Fact \ref{Frobenius}, the application
of Proposition \ref{metacyclic} together
with Fact \ref{Frobenius} will invoke
the Sylow $p$-subgroups of groups of finite
Morley rank for various primes $p$. This
is a subtle issue, because in its full generality,
the Sylow theory in groups of finite Morley
rank is well understood only
for the case $p=2$. This was achieved by Borovik and
Poizat in \cite{BP} where the conjugacy
of the Sylow $2$-subgroups of an arbitrary
group of finite Morley rank was proven.
Moreover, \cite{BP} contained structural
information about these subgroups which
allowed to divide the groups of finite Morley
rank into four mutually disjoint classes:
the groups of {\em odd
type} in which the Sylow $2$-subgroups
are finite extensions of $2$-tori;
the groups of {\em even type}
where the Sylow $2$-subgroups are finite extensions
of infinite unipotent $2$-groups;
the groups of {\em mixed type} where the Sylow
$2$-subgroups are finite extensions of a central
product of a nontrivial unipotent $2$-group and
a nontrivial $2$-torus; the groups of {\em degenerate type}
where the Sylow $2$-subgroups are finite.
An important recent advance
in the theory of groups of finite Morley rank is
the main theorem of \cite{bbc} which states that
a connected group of finite Morley rank of degenerate
type has no involutions. The terminology of types which we
have taken a considerable space to review will be used
below.

Now, we will investigate the effects of Fact \ref{Frobenius}
in our situation.
A finite group is metacyclic if and only if
all of its Sylow $p$-subgroups are cyclic
 \cite[7.6.2 p.~258]{gor}.
So any Frobenius complement is metacyclic if and only if
 it does not contain the quaternion group $Q_8$ of order 8,
which is clearly contained inside
any generalized quaternion group.
In particular, connected
groups of degenerate type have metacyclic Weyl groups
 since they contain no involutions by \cite{bbc}.

\begin{cor}\label{metacyclic_degenerate}
Let $G$ be a connected group of
finite Morley rank of degenerate type,
$T$ any $\pi$-torus in $G$ and
$W=N(T)/C(T)$. If
$C_T(a)$ is finite for all $a \in W^\#$
then $W= N(T)/C(T)$ is metacyclic.
\end{cor}

We begin the next part of our analysis with
 another observation about minimal decent torii, which
 again clarifies the nature of the problem being discussed.
If a minimal decent torus $T$ has a non-cyclic Weyl group,
 then $\pr_p(T) \neq 1$ for any prime $p$
 since $W \hookrightarrow \zz_p^*$ by the above,
and $\zz_p^*$ is cyclic.
It follows that such a $T$ is central in
 any solvable group containing it,
 by the Zilber Field Theorem.

We will need more detailed information about finite subgroups
of algebraic groups over $2$-adic numbers. More precisely,
we will show that (Proposition \ref{noq8ingl2q2})
$\GL_2(\qq_2)$ does not have a subgroup isomorphic
to the group $Q_8$. We imagine that
this is well-known to the specialists. Nevertheless,
we have not been able to find any written reference for this fact
which is provable, as we do below, by elementary means.
In general, there is a conceptual and well developed
approach to such problems concerning subgroups of algebraic
groups over $p$-adic fields (see Proposition 1.12
in \cite{platonrapinnumbers})
which provides the necessary background to analyze
and use the structure
of congruence subgroups (see Chapters 5 and 6 in \cite{analprop}).
Unfortunately, in the case of the prime $2$,
the structure of congruence subgroups is somewhat
exceptional which necessitates additional arguments
for our question. Thus we
have decided to be content with our elementary argument which only
needs basic information about $2$-adic numbers for
which an excellent source is \cite{coursdarithmetique}.

\begin{fact}[{\cite[Chapitres II et III]{coursdarithmetique}}]
\label{hilbertetal}
The following are true in $\qq_2$:
\begin{enumerate}
\item $-1$ is not a square.
\item The equation $x^2+y^2=-1$ does not
have a solution.
\end{enumerate}
\end{fact}
\proof

1. This is part of the Corollaire to Th\'eor\`eme 4
in Chapitre II of \cite{coursdarithmetique}.


2. This uses the Hilbert symbols. Using the notation
in Chapitre III of \cite{coursdarithmetique}, the
statement of 2 can be reformulated as $(-1,-1)=-1$
in $\qq_2$. Th\'eor\`eme 1 of Chapitre III states
that for $(a,b)\in\qq_2^*$
\[(a,b)=(-1)^{\epsilon(u)\epsilon(v)+\alpha\omega(v)+\beta\omega(u)}\]
where $u$ and $v$ are the units in the representation
of $a$ and $b$ respectively used also in point 1.
The function $\epsilon(u)$ is the residue class
of $\frac{u-1}{2}$ modulo $2$. When $a=b=-1$,
this statement becomes
$(-1,-1)=(-1)^{\epsilon(u)^2}$.
But $\epsilon(-1)=-1=1$ modulo $2$.
\qed

\begin{prop}\label{noq8ingl2q2}
The group $\GL_2(\qq_2)$ does not contain
an isomorphic copy of the group of quaternions
$Q_8$.
\end{prop}
\proof Suppose towards a contradiction that
$G\leq\GL_2(\qq_2)$ is isomorphic to $Q_8$.
Let $X\in G$ be a matrix of order $4$. The minimal
polynomial of $X$ is a divisor of $T^4-1$ which
factorizes as $(T^2+1)(T+1)(T-1)$ over $\qq_2$.
By Fact \ref{hilbertetal} 1, these are irreducible
factors. As $|X|=4$, it follows that the minimal
polynomial of $X$ is $T^2+1$. It follows that
$X$ is conjugate to
$\left(\begin{array}{cc}
0&-1\\
1&0
\end{array}\right)$ in $\qq_2$.
Therefore, we may assume that
$X=\left(\begin{array}{cc}
0&-1\\
1&0
\end{array}\right)$ in $\qq_2$.

It follows from the structural properties
of $Q_8$ that there exists
a matrix
$\left(\begin{array}{cc}
a&b\\
c&d
\end{array}\right)$ of order $4$
in $\GL_2(\qq_2)$
that inverts $X$ by conjugation.
By multiplying matrices and using
the fact this second matrix is of
order $4$, we first conclude
that $a=-d$ and $b=c$ and then
$a^2+b^2=-1$. By Lemma \ref{hilbertetal}
this equation does not have a solution
in $\qq_2$, a contradiction.
\qed

We can now constrain the actions upon 2-tori as follows.

\begin{cor}\label{noq8}
Let $G$ be a connected group of finite Morley rank,
 let $T$ be a $\pi$-torus of $G$, and
 consider some $W \leq N(T) / C(T)$.
Suppose that $C_T(a)$ is finite for all $a \in W^\#$.
If either $\pr_2(T)$ is odd or $\pr_2(T) = 2$ then $W$ is metacyclic.
\end{cor}

\proof By Proposition \ref{metacyclic},
$W$ is a Frobenius complement.
Suppose that $W$ is not metacyclic.
Then the discussion following Fact \ref{Frobenius} shows
that there is some $W_0 \cong Q_8$ inside $W$.
We may also take $\pr_2(T) \neq 0$ using our hypotheses.
Set $n= \pr_2(T)$ and $V=\qq_2^n$.
As in the proof of Proposition \ref{metacyclic},
there is a faithful representation of $W_0$ on $V$ by Fact \ref{endptor},
 and this becomes a product of irreducible representations over $\cc$.

As $Q_8$ has 5 conjugacy classes,
 it has 5 irreducible complex representations:
one trivial, three one-dimensional, and the usual two-dimensional one.
Only the two-dimensional one satisfies the hypothesis on
the centralizers of nontrivial elements of $W$.
So our $2$-adic representation becomes
 a product of two-dimensional representations over $\cc$,
and thus $n$ is even. Our assumption then forces $n=2$,
but this is excluded by Proposition \ref{noq8ingl2q2}.
\qed

We now show that all minimal connected simple
groups have metacyclic Weyl groups,
which is the culmination of the preceding preparations:
\\

\begin{lemma}\label{cdtconsolv}
Let $G$ be a connected solvable group
of finite Morley rank and $T$ a
decent torus. Then $N_G(T)$ is connected
and $N_G(T)=C_G(T)$.
\end{lemma}
\proof The first statement
is a corollary of
Facts \ref{classicarter}.1, 2 and 3. Indeed,
let $x\in N(T)$. Then $x$ normalizes
$C^\circ(T)$. By Fact \ref{classicarter}.1,
$C^\circ(T)$ contains a Carter
subgroup $Q$ of $G$. The Frattini argument
applied to $C^\circ(T)$ and
Fact \ref{classicarter}.4 yield
$N(T)=(N(T)\cap N(Q))C^\circ(T)=
QC^\circ(T)=C^\circ(T)$.
\qed

\begin{repstat}{Theorem \ref{cyclic_cor}}
Let $G$ be a minimal connected simple
group of finite Morley rank,
and let $T$ be a decent torus of $G$.
Then the Weyl group
$W= N(T)/C(T)$ is metacyclic.\\
\end{repstat}

In addition to what has been done since the
beginning of this section and the recurrent use
of Theorem \ref{cdtcon}, two substantial facts will
be invoked in the proof. The first one involves
simple groups of odd type:

\begin{fact}[{\cite[High Pr\"ufer Rank Theorem]{burdgsignbalance,bbn,bcj}}]
\label{hipru}
A simple $K^*$-group of finite Morley rank with Pr\"ufer 2-rank at least three
is algebraic.
\end{fact}

\noindent
A {\em $K^*$-group} of finite Morley rank is one in which
every proper, infinite, definable, simple section
is an algebraic group over an algebraically closed
field. Note that a minimal counterexample to the
Cherlin-Zil'ber conjecture is a $K^*$-group.
The second fact involves simple groups of even type:

\begin{fact}[{\cite{abcbook}}]\label{eventheorem}
A simple group of finite Morley rank which contains
an infinite elementary abelian $2$-subgroup is
a linear algebraic group over an algebraically closed
field of characteristic $2$.
\end{fact}

We have accumulated
sufficient arsenal to carry out the proof
of Theorem \ref{cyclic_cor}. It should be reminded
that our key result in this paper, namely
Theorem \ref{cdtcon}, will keep playing its important
but quiet role since it is the reason why
the centralizers of all decent or $\pi$-tori
are connected although this crucial observation
is not acknowledged every time it is used.\\

\proofof{Theorem \ref{cyclic_cor}}
If $G$ has even type, then $G \cong SL_2(k)$
 where $k$ is an algebraically closed field
of characteristic 2 by Fact \ref{eventheorem}.
So we may assume that $G$ is of odd or degenerate type.

The main part of the argument will not distinguish
between odd and degenerate types.
Let $S$ be a maximal decent torus of $C(T)$, and hence of $G$.
By a Frattini argument which
uses Fact \ref{dtgen}.2, $N(T) = N_{N(T)}(S) C(T)$.
Clearly $C(S) \leq C(T)$ since $T \leq S$.
So $W$ is a section of $N(S) / C(S)$.
We may thus assume that $T$ is a maximal decent torus of $G$,
 in particular $\pr_2(G) = \pr_2(T)$.

For any decent torus $R \leq T$,  
 $C(R)$ is connected by Theorem \ref{cdtcon},
  and hence solvable.
So $N_{C(R)}(T) = C_{C(R)}(T)$ by Lemma \ref{cdtconsolv}.
Thus any $a \in W \cap C(R)$ lies inside $W \cap C(T)$ too,
 but this second group is trivial.
Hence $C_T(w)$ has finite torsion for any
$w\in W$.

By Corollary \ref{metacyclic_degenerate}, we are done if $G$ is
of degenerate type. Thus,
by Corollary \ref{noq8}, either $W$ is metacyclic
 or else $\pr_2(T) \geq 4$.
But $\pr_2(G) \leq 2$ by the High Pr\"ufer Rank Theorem,
 since $G$ is minimal connected simple.
\qed

We observe that here one may replace the use of
the High Pr\"ufer Rank Theorem,
as well as part of the argument,
by an application of Adrien Deloro's thesis
\cite{delorothesis}.
In a similar vein, one does not need the full force
of Fact \ref{eventheorem} but the analysis
of minimal simple connected groups of finite Morley
rank, which was nevertheless the most demanding
part of the classification of simple groups
of even type.
We also observe that Theorem \ref{cyclic_cor}
should hold under Hypothesis \ref{ngroup},
assuming that Deloro and Jaligot's program
(see the beginning of
Section \ref{minbackg}) is successful in that context.

\section{Automorphisms of minimal simple groups}
\label{autosection}

Unlike the High Pr\"ufer Rank Theorem (Fact \ref{hipru}),
the classification simple groups of finite Morley
rank of even type (Fact \ref{eventheorem}) is independent
of any inductive hypothesis. Its achievement has recently
motivated a considerable body of work in the analysis of
simple groups of odd type around a similar strategy
which consists of carrying out an inductive proof with
weaker inductive assumptions.
The second author and Cherlin have taken
 the initiative to move in this direction.
The Deloro-Jaligot enterprise mentioned in Section \ref{minbackg}
 is part of the same picture.

The analysis of simple groups of odd type is more problematic
than that of simple groups of even type in that
in the latter case the initial step in the induction
was rather simple due to the very restricted behavior
of unipotent 2-groups of automorphisms of minimal
simple groups.  In the case of simple groups of odd type
the action of a 2-torus over a minimal connected
simple group is far from being understood,
and any new technique or tool is worth considering.
In this vein, we analize automorphisms of those connected
minimal simple groups of finite Morley rank which
most resemble algebraic groups.
Our viewpoint here is that Theorem \ref{cdtcon} forbids
outer toral automorphisms in the presence of a Weyl group.

\smallskip

As will be rapidly
verified in Lemma \ref{verifyngroup} below,
Hypothesis \ref{ngroup} offers the right setting for the
intended analysis.
The following conditions describe the particular
properties which lead to Hypothesis \ref{ngroup}.

\begin{hypothesis}\label{hypothesis:global}
Let $\hat{G} = G\sd A$ be a group of finite Morley rank where
\begin{itemize}
\item G is a definable minimal connected simple subgroup
normal in $\hat G$,
\item A is a connected definable abelian subgroup, and
\item $\hat{G}$ is centerless.
\end{itemize}
\end{hypothesis}

Before stating our main theorem, we verify that
the group $\hat G$ in Hypothesis \ref{hypothesis:global}
indeed implies Hypothesis
\ref{ngroup}.

\begin{lemma}\label{verifyngroup}
Let $\hat G$, $G$ and $A$ be as in the
statement of Hypothesis \ref{hypothesis:global}.
Then $\hat G$ satisfies Hypothesis \ref{ngroup}.
\end{lemma}
\proof Clearly $\hat G$ is not solvable. Let $U$
be a nontrivial, definable, abelian
subgroup of $\hat G$ such that $X=N_{\hat G}^\circ(U)$
is not solvable.
Since $A$ is abelian, $X\cap G\neq1$ and
$X\cap G$ is not solvable.
Since $G$ is minimal, either $X\cap G$ is
finite or $G\leq X$.

If $X\cap G$ is finite, then since
$X'\leq X\cap G$, $X'$ is connected
and finite. Thus $X'$ is trivial which contradicts
that $X$ is not solvable.

If $G\leq X$ then
$[G,U]=1$. This contradicts that $Z(\hat G)=1$.
\qed

The following is the main technical step towards
the proof of Theorem \ref{autthm_torii}:

\begin{prop}\label{genaut}
Let $\hat G = G \sd A$ be a group of finite Morley
rank where both $G$ and $A$ are definable connected
subgroups, and $G$ is normal in $\hat G$
while $A$ is abelian.
Let $Q$ be a definable connected nontrivial
subgroup of $C_G(A)$.
Suppose that $\hat G/G$ has nontrivial
divisible torsion. Then
every coset
of $N_G(Q)/N_G^\circ(Q)$ is represented
by an element $w$ to which
is associated a $p$-torus $T_w$ of $\hat G$
such that $\langle w, Q \rangle\leq C_G(T_w)$.
\end{prop}
\proof
By hypothesis, $\hat G/G$ has nontrivial
divisible torsion.
By Fact \ref{nodeftorsion}, $A$ has
nontrivial divisible torsion which is not
contained in $G$.

Since by assumption $[A,Q]=1$,
$C^\circ_{\hat G}(Q)$ has a nontrivial $p$-torus $T$
for some prime $p$,
which we take maximal in $C_{\hat G}(Q)$. Furthermore,
we may assume that $A = d(A \cap T)$.
In particular, $T$ and $A$ commute.
So, by a Frattini argument (Fact \ref{dtgen}.2) applied
to the conjugates of $T$ in $N_{\hat G}(Q)$,
\[N_{\hat G}(Q) = N_{N_{\hat G}(Q)}(T) C^\circ_{\hat G}(Q).
\ \ \ (1) \]
On the other hand, since $\hat G=GA$ and
$A\leq C_{\hat G}(Q)$ by assumption,
we have
\[N_{\hat G}(Q)=N_{\hat G}(Q)\cap GA=
N_G(Q)A=N_G(Q)C_{\hat G}^\circ (Q).\ \ \ (2)\]

It follows from (1) and (2)
and the fact that $[A,T]=1$ that
for any $x \in N_{\hat G}(T)\cap N_{\hat G}(Q)$,
 there exists $y\in G$ such that
$x A=yA$ where
$y\in N_G(Q)\cap N_{\hat G}(T)$.
Thus $y\in N_G(Q)\cap N_G(T)$.
It follows from this
that
\[ N_{N_{\hat G}(Q)}(T) C^\circ_{\hat G}(Q)=
N_{N_G(Q)}(T) C^\circ_{\hat G}(Q).\]
On the other hand,
since we have already shown
that
\[ N_{N_{\hat G}(Q)}(T) C^\circ_{\hat G}(Q)
=N_G(Q) A = N_G(Q) C^\circ_{\hat G}(Q),\]
we conclude
$N_{N_G(Q)}(T) C^\circ_G(Q) = N_G(Q)$.



Our initial setup yields $[TG:G]=\infty$.
Let $w$ be a representative of a coset
of $N_G^\circ(Q)$ in $N_G(Q)$.
Since $[w,T]\leq T\cap G<T$,
we thus conclude that $C_T(w)$ is infinite.
We let $T_w=C_T^\circ(w)$ 
and $H=C_G(T_w)$.
We have $w\in H$ and $Q \leq C^\circ_G(T) \le H$.
\qed

We will apply Proposition \ref{genaut} in the
minimal context of Hypothesis \ref{hypothesis:global}
to obtain the following corollary from which
the main theorem of this section will immediately
follow.

\begin{cor}\label{autthm}
Let $\hat G$, $G$ and $A$ be as in
Hypothesis \ref{hypothesis:global}.
Let $Q$ be a definable connected nontrivial
subgroup of $C_G(A)$ such that
 for any definable $H \leq G$ containing $Q$,
 the group $N^\circ_H(Q)$ is almost self-normalizing.
Suppose that $\hat G/G$ has nontrivial
divisible torsion.
Then $N_G(Q)$ is connected. 
In particular,
 if $Q$ lies inside a Carter subgroup $Q_1$ of $G$,
 then $Q_1$ is self-normalizing. 
\end{cor}

In the proof of Corollary \ref{autthm}, we will need a
consequence of Fr\'econ's
analysis of abnormal subgroups of connected
solvable groups of finite Morley rank.

\begin{fact}[{\cite[Th\'eor\`eme 1.2]{frec1}}]
\label{freconabnormal}
Let $G$ be a connected solvable group of finite
Morley rank and $H$ a definable subgroup of $G$.
Then the following are equivalent:
\begin{itemize}
\item $H$ is abnormal
\item $H$ is definable and connected, and there
exist $n\in\nn^*$ and a decreasing sequence
of subgroups $(H_i)_{0\leq i\leq n}$ such that
$H_0=G$, $H_n=H$ and for $i>0$ each $H_i$ is
a proper definable connected maximal subgroup
of $H_{i-1}$ not normal in $H_{i-1}$. 
\end{itemize}
\end{fact}
A subgroup $H$ of an arbitrary group $G$ is
said to be {\em abnormal} if $g\in\langle H,H^g\rangle$
for every $g\in G$.

\begin{lemma}\label{selfnormsolv}
Let $G$ be a connected solvable group
of finite Morley rank and $H\leq G$ be
a definable connected subgroup of $G$
such that $[N_G(H):H]<\infty$. Then
$N_G(H)=H$.
\end{lemma}
\proof If $H=G$ then there is nothing
to do. If $H<G$ then a sequence as
in Fact \ref{freconabnormal} can be constructed
so that the smallest term is $H$.
Thus $H$ is abnormal and hence
self-normalizing.
\qed

\proofof{Corollary \ref{autthm}}
The hypotheses of Proposition \ref{genaut}
are clearly satisfied. Thus, if
$w\in N_G(Q)\setminus N_G^\circ(Q)$, then
there exists a $p$-torus $T_w$ in $\hat G$
which centralizes $w$ and $Q$.
By Theorem \ref{cdtcon}, $C_{\hat G}(T_w)$ is connected.
Since $\hat G=G\sd A$, it follows that
$H$ is connected.
We have $w\in H$ and $Q \leq C^\circ_G(T) \le H$.
By our hypothesis on $Q$,
$N_H^\circ(Q)$ is almost self-normalizing in $H$,
and hence self-normalizing by Lemma \ref{selfnormsolv}.
So $w \in N^\circ_H(Q) \leq N^\circ_G(Q)$, a contradiction.
\qed

In particular,
 taking $Q$ to be a decent torus in a minimal
simple group $G$
 yields the target theorem of this section:\\

\begin{repstat}{Theorem \ref{autthm_torii}}
Let $\hat G$, $G$ and $A$ be as in Hypothesis
\ref{hypothesis:global}. Then
the following two properties are contradictory.
\begin{itemize}
\item  $G$ has a divisible abelian torsion subgroup $Q$
 with a non-trivial Weyl group $W= N_G(Q)/C_G(Q)$.
\item  $\hat{G} / G$ has nontrivial divisible torsion.
\end{itemize}
\end{repstat}

\bibliography{../../Texdestek/bibl}

\begin{thebibliography}{10}

\bibitem{baudgroup}
A.~Baudisch.
\newblock A new uncountably categorical group.
\newblock {\em Trans. Amer. Math. Soc.}, 348:3889--3940, 1996.

\bibitem{ber1}
A.~Berkman.
\newblock The classical involution theorem for groups of finite morley rank.
\newblock {\em J. Algebra}, 243:361--384, 2001.

\bibitem{berbor01}
A.~Berkman and A.~Borovik.
\newblock A generic identification theorem for groups of finite {M}orley rank.
\newblock {\em J. London Math. Soc.}, 69 (1):14--26, 2004.

\bibitem{bbc}
A.~Borovik, J.~Burdges, and G.~Cherlin.
\newblock Involutions in groups of morley rank of degenerate type.
\newblock submitted.

\bibitem{bbn}
A.~Borovik, J.~Burdges, and A.~Nesin.
\newblock Uniqueness cases in odd type groups of finite {M}orley rank.
\newblock {\em J. London Math. Soc.}
\newblock to appear.

\bibitem{bnbook}
A.~V. Borovik and A.~Nesin.
\newblock {\em Groups of Finite Morley Rank}.
\newblock Oxford University Press, 1994.

\bibitem{BP}
A.~V. Borovik and B.~Poizat.
\newblock Tores et $p$-groupes.
\newblock {\em J. Symbolic Logic}, 55:478--490, 1990.

\bibitem{burbender}
J.~Burdges.
\newblock The {B}ender method in groups of finite {M}orley rank.
\newblock Submitted.

\bibitem{burdgsignbalance}
J.~Burdges.
\newblock Signalizers and balance in groups of finite {M}orley rank.
\newblock submitted.

\bibitem{bcj}
J.~Burdges and G.~Cherlin a~E.~Jaligot.
\newblock Minimal connected simple groups of finite {M}orley rank with strongly
  embedded subgroups.
\newblock {\em J. Algebra}.
\newblock to appear.

\bibitem{bcsstorsion}
J.~Burdges and G.~Cherlin.
\newblock Semisimple torsion in groups of finite morley rank.
\newblock in preparation.

\bibitem{goodtori}
G.~Cherlin.
\newblock Good tori in groups of finite morley rank.
\newblock {\em J. Group Theory}, 8:613--621, 2005.

\bibitem{delorothesis}
A.~Deloro.
\newblock {\em Groupes simples connexes minimaux de type impair}.
\newblock PhD thesis, Universit\'e Paris 7 - Denis Diderot, 2007.

\bibitem{analprop}
J.~D. Dixon, M.~P.~F. du~Sautoy, A.~Mann, and D.~Segal.
\newblock {\em Analytic Pro-p Groups}.
\newblock Cambridge University Press, second edition, 1999.

\bibitem{frecpseudotori}
O.~Fr\'econ.
\newblock Pseudo-tori and subtame groups of finite {M}orley rank.
\newblock in preparation.

\bibitem{frec1}
O.~Fr\'econ.
\newblock Sous-groupes anormaux dans les groupes de rang de {M}orley fini
  r\'esolubles.
\newblock {\em J. Algebra}, 229:118--152, 2000.

\bibitem{frecjalcar}
O.~Fr\'econ and E.~Jaligot.
\newblock The existence of {C}arter subgroups in groups of finite {M}orley
  rank.
\newblock {\em J. Group Theory}, 8:623--633, 2006.

\bibitem{frecjalnewt}
O~Fr\'econ and E.~Jaligot.
\newblock Conjugacy in groups of finite {M}orley rank.
\newblock Cambridge University Press, 2007.

\bibitem{FuHa}
William Fulton and Joe Harris.
\newblock {\em Representation Theory}.
\newblock Number 129 in Graduate Texts in Mathematics. Springer Verlag, 1991.

\bibitem{gor}
D.~Gorenstein.
\newblock {\em Finite Groups}.
\newblock Chelsea Publishing Company, New York, 1980.

\bibitem{hu}
J.~E. Humphreys.
\newblock {\em Linear Algebraic Groups}.
\newblock Springer Verlag, Berlin-New York, second edition, 1981.

\bibitem{ba1}
N.~Jacobson.
\newblock {\em Basic Algebra I}.
\newblock W. H. Freeman and Company, 1989.

\bibitem{ba2}
N.~Jacobson.
\newblock {\em Basic Algebra II}.
\newblock W. H. Freeman and Company, 1989.

\bibitem{jacart}
E.~Jaligot.
\newblock Generix never gives up.
\newblock {\em J. Symbolic Logic}, 71:599--610, 2006.

\bibitem{abcbook}
T.~Alt\i nel, A.~Borovik, and G.~Cherlin.
\newblock {\em Simple Groups of Finite Morley Rank}.
\newblock American Mathematical Society, expected 2007.

\bibitem{N1}
A.~Nesin.
\newblock On solvable groups of finite {M}orley rank.
\newblock {\em Trans. Amer. Math. Soc.}, 321:659--690, 1990.

\bibitem{N3}
A.~Nesin.
\newblock Poly-separated and $\omega$-stable nilpotent groups.
\newblock {\em J. Symbolic Logic}, 56:694--699, 1991.

\bibitem{platonrapinnumbers}
Vladimir Platonov and Andrei Rapinchuk.
\newblock {\em Algebraic Numbers and Number Theory}.
\newblock Academic Press, 1994.

\bibitem{cours}
B.~Poizat.
\newblock {\em Cours de Th\'eorie des Mod\`eles}.
\newblock Nur Al-Mantiq Wal-Ma'rifah, Villeurbanne, France, 1985.

\bibitem{grst}
B.~Poizat.
\newblock {\em Groupes Stables}.
\newblock Nur Al-Mantiq Wal-Ma'rifah, Villeurbanne, France, 1987.

\bibitem{coursdarithmetique}
J.-P. Serre.
\newblock {\em Cours d'arithm\'etique}.
\newblock Presses Universitaires de France, Paris, 1970.

\bibitem{carwag}
F.~Wagner.
\newblock Nilpotent complements and {C}arter subgroups in stable {$\frak
  R$}-groups.
\newblock {\em Arch. Math. Logic}, 33:23--34, 1994.

\bibitem{stgr}
F.~Wagner.
\newblock {\em Stable Groups}.
\newblock Number 240 in London Mathematical Society Lecture Note Series.
  Cambridge University Press, Cambridge, UK, 1997.

\bibitem{wagbadp}
F.~Wagner.
\newblock Fields of finite {M}orley rank.
\newblock {\em J. Symbolic Logic}, 66:703--706, 2001.

\end{thebibliography}
\bibliographystyle{plain}

\end{document}